\documentclass[11pt]{article}
\usepackage{amssymb}
\usepackage{amsmath}

\newtheorem{theo}{Theorem}[section]

\newtheorem{prop}[theo]{Proposition}
\newtheorem{lemm}[theo]{Lemma}
\newtheorem{coro}[theo]{Corollary}
\newtheorem{rema}[theo]{Remark}
\newtheorem{Defi}[theo]{Definition}
\newtheorem{ex}[theo]{Example}

\newcommand{\cqfd}
{%
\mbox{}%
\nolinebreak%
\hfill%
\rule{2mm}{2mm}%
\medbreak%
\par%
}
\newfont{\gothic}{eufb10}
\date{\empty}
\begin{document}
\title{Hodge structures on cohomology algebras and geometry}
\author{Claire Voisin\\ Institut de math{\'e}matiques de Jussieu, CNRS,UMR
7586} \maketitle \setcounter{section}{-1}
\section{Introduction}
\setcounter{equation}{0}
It is well-known (see eg \cite{voisinbook}) that the topology of a compact K\"ahler manifold
$X$ is strongly restricted
by Hodge theory. In fact, Hodge theory provides two sets of data on the cohomology of
a compact K\"ahler manifold. The first data are the Hodge decompositions
on the cohomology spaces $H^k(X,\mathbb{C})$
(see (\ref{hodgedecomp}) where $V=H^k(X,\mathbb{Q})$); they
 depend only on the complex structure.

 The second data, known as the Lefschetz isomorphism and the Lefschetz decomposition
 on cohomology (see (\ref{lefdecomp}) with $A^k_\mathbb{R}=H^k(X,\mathbb{R})$) depend only on the choice
 of a K\"ahler class, but remain satisfied by any symplectic class close to a
 K\"ahler class.

 Both are combined to give the so-called  Lefschetz bilinear relations, which
 lead for example to the Hodge index theorem (cf \cite{voisinbook}, 6.3.2) which computes the signature
 of the intersection form on the middle cohomology of an even dimensional compact K\"ahler
 manifold as an alternate sum of its Hodge numbers.

 If we want to extract topological restrictions using these informations, we are faced to
 the following problem: neither the complex structure (or even the Hodge numbers), nor the (deformation
 class) of the symplectic structure (or even the symplectic class) are topological.

 For this reason, only a very small number of purely topological restrictions
 have been extracted so far from these data. (Note however that the formality theorem
 \cite{dgms},
 which uses more than the data above, is a topological statement. Similarly, non-abelian Hodge theory has
 provided strong restrictions on $\pi_1(X)$ (cf \cite{burgeramoros}).)
The classically known restrictions are the following:

\begin{enumerate}
\item Due to the Hodge decomposition and the Hodge symmetry (\ref{hodgesymmetry}),
the odd Betti numbers $b_{2i+1}(X)$ have to be even.
\item Due to the Lefschetz property, the  even Betti numbers
$b_{2i}(X)$ are increasing in the range $2i\leq n=dim_\mathbb{C}X$ and similarly the
odd Betti numbers $b_{2i+1}(X)$ are increasing in the range
$2i+1\leq n=dim_\mathbb{C}X$.
\end{enumerate}
(Note that, as the dimension of the manifold
is even,  the Lefschetz property also implies the condition of evenness of odd Betti numbers.)

The purpose of this paper is to extract from the Hodge decomposition and the Lefschetz property
a number of purely topological restrictions on the {\it cohomology algebra}
of a compact K\"ahler manifold. We will show that these restrictions are effective
even in the category of
compact symplectic manifold satisfying the Lefschetz property.

In \cite{thurston}, \cite{mcduff}, examples of compact symplectic manifolds
which are topologically non K\"ahler were constructed. The examples did not have
their odd Betti numbers $b_{2i+1}(X)$ even. In \cite{gompf} and \cite{bouyakoub}, one can find
 examples of compact symplectic manifolds
whose odd  Betti numbers are even, but for which  the Lefschetz property is satisfied by no degree
$2$ cohomology class.

Here, many of  the examples we construct satisfy the Lefschetz property. Furthermore, they are all
built starting from compact K\"ahler manifolds, and considering either
symplectic blow-up of them in a ``wrong'' symplectic embedding, or
complex projective bundles on them. For all of them we conclude that
their rational (and even sometimes real) cohomology algebra does not satisfy the restrictions
imposed to the cohomology algebra of a compact K\"ahler manifold.

As in \cite{voisin}, the key point here is the observation that the Hodge decomposition
on the cohomology of a compact K\"ahler manifold is compatible with the cup-product
(see (\ref{compatcond})). This leads to a number of algebraic restrictions on the cohomology algebra
of $X$. Note that what we provide here is only a sample of them, where we tried
to separate restrictions of three kinds:

\begin{enumerate}
\item Restrictions on the real cohomology algebra coming from the Hodge decomposition.
\item In the spirit of \cite{voisin}, more subtle restrictions on the rational cohomology algebra
coming from the Hodge decomposition.
\item Restrictions coming from the polarization on the Hodge structure.
\end{enumerate}

These restrictions, together with examples showing that they are all
 effective even in the  symplectic category, are described in section \ref{sec3}.
 In section \ref{sec2},
 we give two stability results concerning cohomology algebras endowed with a polarized
 Hodge structure. The first one (Theorem \ref{product}) concerns tensor products of cohomology algebras (corresponding
 to taking products of manifolds):
 \begin{theo} Assume there is a polarized Hodge structure on a
 (rational or real) cohomology algebra
$M$, and assume that
$$M\cong A\otimes B,$$
where $A$ and $B$ are (rational or real) cohomology algebras. If either
$A^1=0$ or $B^1=0$, then $A$ and $B$ are of even dimension
and there are polarized Hodge structures on
$A$ and $B$, inducing that of $M$.
 \end{theo}
 The other one  (Theorem \ref{theoprojbun}) concerns projective bundles :

 \begin{theo} Let $X$ be a compact connected smooth oriented  manifold,  and let
$E$ be a complex vector bundle on $X$ with trivial determinant. Assume that the cohomology
of $X$ is generated in degrees $1$ and $2$. Then if the cohomology algebra
$H^*(\mathbb{P}(E),\mathbb{Q})$ carries a  Hodge structure, the cohomology algebra
$H^*(X,\mathbb{Q})\subset H^*(\mathbb{P}(E),\mathbb{Q})$ has an induced Hodge  structure, for which
the Chern classes $c_i(E)$ are Hodge classes. A similar result holds with
$H^*(\mathbb{P}(E),\mathbb{Q}), H^*(X,\mathbb{Q})$
replaced by $H^*(\mathbb{P}(E),\mathbb{R}),\,H^*(X,\mathbb{R})$.
 \end{theo}
These results are used in section \ref{sec3} to construct
compact symplectic manifolds with non-K\"ahler rational cohomology algebras, but satisfying
the Lefschetz property. In one of our examples, the criterion we use, namely
the existence of a Hodge structure on the cohomology algebra, needs the {\it rational}
cohomology algebra, while on  the first two examples, the real cohomology algebra
suffices to exclude the existence of a Hodge structure.

In the last section, we show that there are in fact supplementary
constraints on the cohomology
algebra of a compact K\"ahler manifold.
Indeed, we finally construct an example of a compact symplectic manifold which has
the {\it real} cohomology algebra of a K\"ahler manifold,
and whose cohomology algebra carries a
rational Hodge structure, but whose rational cohomology algebra is not
the rational cohomology algebra of a compact K\"ahler manifold.
Thus all the restrictions we met before are satisfied, but still we find supplementary
restrictions of a more subtle nature, related to the result in \cite{voisinimrn}:
some Hodge classes on certain compact K\"ahler manifolds cannot be constructed as
Chern classes of holomorphic vector bundles or even analytic coherent sheaves.

\section{Real and rational Hodge structures on cohomology algebras}
Let us recall the notion of a Hodge structure \cite{voisinbook}, 7.1.1.
\begin{Defi} A rational (resp. real) Hodge structure of weight $k$ is
a rational (resp. real) vector space $V$, together with a decomposition into a direct sum of complex
vector  subspaces
\begin{eqnarray} \label{hodgedecomp}
V_\mathbb{C}:=V\otimes\mathbb{C}=\oplus_{p+q=k}V^{p,q},\end{eqnarray}
satisfying the Hodge symmetry condition
\begin{eqnarray}\label{hodgesymmetry}\overline{V^{p,q}}=V^{q,p}.\end{eqnarray}
\end{Defi}
Here the tensor product is taken over $\mathbb{Q}$ in the rational case, and over $\mathbb{R}$
in the second case.
\begin{rema}\label{remapaire}{\rm A very classical observation is the fact that
for odd $k$, the existence of a Hodge structure of weight $k$ on $V$ forces the dimension of
$V$ to be even, by the Hodge symmetry (\ref{hodgesymmetry}).}
\end{rema}
As is well known, the data of a Hodge structure of weight $k$ on $V$
 is equivalent to the data of an action
of $\mathbb{C}^*$ on $V_\mathbb{R}$, satisfying the property that
$\lambda\in \mathbb{R}^*$ acts by multiplication by $\lambda^k$. Indeed, we let
$z\in \mathbb{C}^*$ act on $V_\mathbb{C}$ by multiplication by
$z^p\overline{z}^q$ on $V^{p,q}$ and the Hodge symmetry (\ref{hodgesymmetry})
implies that this action leaves $V_\mathbb{R}:=V\otimes\mathbb{R}$ stable.

In this paper, we will consider what we will  call
rational (resp. real) cohomology algebras, that is finite dimensional
graded associative  $\mathbb{Q}$-algebras  (resp. $\mathbb{R}$-algebras)
with unit, satisfying the conditions
making them good candidates to be cohomology algebras of connected compact oriented manifolds:
\begin{enumerate}
\item The product is graded commutative.
\item \label{cond212mars}$A^i=0$ for $i<0$ and the  term $A^0$ is generated over $\mathbb{Q}$ (resp. $\mathbb{R}$)
by $1_A$.
\item \label{cond312mars} For a certain integer $m$ that we will call the dimension of $A$, we have
$$A^m\cong \mathbb{Q},\,\,{\rm resp.}\,\,A^m\cong \mathbb{R},$$
and the pairing
\begin{eqnarray}\label{*pairing}A^k\otimes A^{m-k}\rightarrow A^m
\end{eqnarray}
is perfect, for any $k$ (so, in particular $A^k=0$ for $k>m$).
\end{enumerate}

Note that
condition \ref{cond212mars} reflects  the connectivity of $X$, when  $A=H^*(X,\mathbb{Q})$,
where $X$ is a topological space, and that the dimension of $A$
is   the dimension of $X$ if  $X$ is a compact oriented  manifold.

We will denote $\cup$ the product on $A$, whether  $A$ comes from geometry or not.
\begin{Defi} A  Hodge structure on a cohomology algebra $A$ is the data of a Hodge structure
of weight $k$ on each graded piece $A^k$, satisfying the following compatibility
property with the product of $A$:

\begin{eqnarray}\label{compatcond} A^{p,q}\cup A^{p',q'}\subset A^{p+p',q+q'}.
\end{eqnarray}
\end{Defi}
\begin{rema}{\rm If there exists a Hodge structure on a cohomology algebra $A$, then the dimension of
$A$ is even. Indeed, the top degree cohomology of $A$ is endowed with a Hodge structure
of weight $l=dim\,A$. On the other hand, it is of rank $1$ by condition \ref{cond312mars}.
By remark \ref{remapaire}, it follows that $l$ is even.}
\end{rema}

The main result of Hodge theory applied to K\"ahler geometry
 is the following (\cite{voisinbook}, 6.1.3):
\begin{theo} Let $X$ be a compact K\"ahler manifold.
Then the Hodge decomposition on each
$H^k(X,\mathbb{Q})$ equips the cohomology algebra $H^*(X,\mathbb{Q})$   with a Hodge structure.
\end{theo}
Our goal in this paper is to explore the {\it topological} restrictions deduced from
the existence of a Hodge structure on $H^*(X,\mathbb{Q})$ or
$H^*(X,\mathbb{R})$. We will show that these restrictions are
effective even in the category of
compact symplectic manifolds satisfying the Lefschetz property (see next section).
It turns out that stronger restrictions are obtained using the notion of polarized Hodge structure.
\subsection{Polarizations}
Another important property satisfied by the cohomology algebra of a compact K\"ahler
manifold $X$ is the fact that the Hodge structure on it can be polarized using a
K\"ahler class $\omega\in H^{1,1}(X)_\mathbb{R}$.
Let us define a polarization on a cohomology algebra endowed with a Hodge structure
$A$. The class of the polarization should be an element
$$\omega\in A^{1,1}_\mathbb{R}:=A^{1,1}\cap A_\mathbb{R}^2.$$

This element, seen as an element of $A^2$, should satisfy the Lefschetz property: let $dim\,A=2n$.
Then for any integer $k$, $0\leq k\leq n$, the cup-product by $\omega^{n-k}$
$$\cup\omega^{n-k}: A^k_\mathbb{R}\rightarrow A^{2n-k}_\mathbb{R}$$
should be an isomorphism. Note that both sides have the same dimension by
the duality (\ref{*pairing}).
\begin{rema} {\rm The morphism $\cup\omega^{n-k}$ is anti-self-adjoint for $k$ odd, and
self-adjoint for $k$ even,
with respect to the duality $A^k_\mathbb{R}\cong (A^{2n-k}_\mathbb{R})^*$. Thus, the existence of a
degree $2$
class $\omega$ satisfying the Lefschetz property above implies that the dimension of
$A^k_\mathbb{R}$ is even for odd $k$.}
\end{rema}
The Lefschetz property implies the Lefschetz decomposition (\ref{lefdecomp}) below
 (cf \cite{voisinbook}, 6.2.3):
Let $\omega$ be a class satisfying the Lefschetz property above, and define for
$k\leq n$ the primitive part
$A^k_{\mathbb{R},prim}$ by
$$A^k_{\mathbb{R},prim}:=Ker\,(\cup\omega^{n+1-k}:
 A^k_\mathbb{R}\rightarrow A^{2n-k+2}_\mathbb{R}).$$
 Then we have for $k\leq n$
\begin{eqnarray}\label{lefdecomp}
\bigoplus_{k-2i\geq0}A^{k-2i}_{\mathbb{R},prim}\stackrel{\sum_i \cup\omega^i}{\cong}A^k_\mathbb{R}.
\end{eqnarray}
Observe that the $A^{k}_{\mathbb{R},prim}$ are real sub-Hodge structures of
$A^{k}_\mathbb{R}$, which means that the corresponding complex vector spaces
$$A^{k}_{\mathbb{C},prim}\subset A^{k}_{\mathbb{C}}$$
are stable under the Hodge decomposition. This follows from condition (\ref{compatcond}) above, and from
the fact that $\omega$ is of type $(1,1)$.
If furthermore $A$ is a rational cohomology  algebra and we can choose $\omega$ to be rational, then
$A^{k}_{\mathbb{R},prim}$ is in fact defined over $\mathbb{Q}$
and provides a rational sub-Hodge structure of $A^{k}$.

To conclude, let us mention the  Riemann bilinear relations, which will play a role in section
\ref{sectionpol}.
$A$ being as above a cohomology algebra of dimension $2n$ endowed
with a Hodge structure and a class
 $\omega$ of type $(1,1)$ satisfying  Lefschetz property, we can construct
for each $k\leq n$, using  the duality on $A$ given by
the generator $\omega^n$ of $A^{2n}_\mathbb{R}$ , a non degenerate intersection pairing
$q_\omega$ on $A^{k}_\mathbb{R}$, which is symmetric if $k$ is even and alternate if $k$ is
odd, defined by:
$$q_\omega(\alpha,\beta)=\omega^{n-k}\cup \alpha\cup\beta\in A^{2n}_\mathbb{R}\cong \mathbb{R}.$$
(The last isomorphism here is defined up to a multiplicative coefficient).

It follows easily from the definition of the primitive parts $A^i_{\mathbb{R},prim}\subset
A^i_\mathbb{R}$ that
the Lefschetz decomposition is orthogonal for the pairing $q_\omega$.
Let us now introduce the Hermitian pairing on $A^k_\mathbb{C}$:
$$h_\omega(\alpha,\beta)=\iota^kq_\omega(\alpha,\overline{\beta}).$$
For bidegree reasons, using the condition (\ref{compatcond}), we also find that the Hodge decomposition
is orthogonal with respect to the
pairing $h_\omega$.
\begin{eqnarray}\label{firstriemann} h_\omega(\alpha,\beta)=0,\,\alpha\in A^{p,q}_\mathbb{C},\,
\beta\in A^{p',q'}_\mathbb{C},\,(p,q)\not=(p',q').
\end{eqnarray}
Finally the second Riemann bilinear relations are  restrictions on the signs
of the Hermitian pairing $h_\omega$
restricted to the part
$$A^{p,q}_{\mathbb{C},prim}\subset A^{p,q}_\mathbb{C},\,p+q=k$$
defined as $A^{p,q}_{\mathbb{C},prim}=A^{p,q}_\mathbb{C}\cap A^{k}_{\mathbb{C},prim}$.
If $A$ is the cohomology algebra of a compact K\"ahler manifold,
these restrictions are described in
the following theorem.
\begin{theo} \label{secondriemann} Let $X$ be a compact K\"ahler
manifold with K\"ahler class
$\omega$.   Then the Hermitian form $h_\omega$ is definite of sign
$(-1)^{\frac{k(k-1)}{2}}\iota^{p-q-k}$ on the component
$$\omega^r\cup H^{p,q}(X,\mathbb{C})_{prim},\,2r+p+q=k$$ of $H^{k}(X,\mathbb{C})$.
\end{theo}
A polarized  cohomology algebra $(A,\omega)$ should satisfy
also these sign restrictions.

At this point, we find another constraint
on the Betti numbers of a compact K\"ahler manifold, which shows that we cannot consider
separately the odd and even Betti numbers to address the question asked by Simpson in
\cite{simpson}, namely what can be the Betti or Hodge  numbers of compact K\"ahler manifolds.
Indeed, we have the following lemma:
\begin{lemm}\label{simpsonproblem} Let $M$ be a cohomology algebra which carries a polarized
Hodge structure. Let $2n=dim\,M$,  and let $i$ be an integer such that
$$4i+2\leq n.$$
Then if $M^{2i+1}\not=0$, we have $rk\,M^{4i+2}\geq2$ and more precisely
$rk\,M^{2i+1,2i+1}\geq2$.
\end{lemm}
{\bf Proof.} Indeed, note  that some $M^{p,q}\not=\{0\}$ for some $p>q,\,p+q=2i+1$. Then
there is an
$\alpha\in M_\mathbb{C}^{p,q}$
such that $\alpha\cup\overline{\alpha}\not=0$, otherwise this contradicts the fact that for a polarisation
given by $\omega\in M^{1,1}$,
the Hermitian form
$$h_\omega(\alpha,\overline{\alpha})=\iota\omega^{n-2i-1}\cup\alpha\cup\overline{\alpha}$$
is non degenerate on $M_\mathbb{C}^{p,q}$ by the Lefschetz decomposition
(\ref{lefdecomp}) and the
 Riemann bilinear relations (\ref{secondriemann}).

Hence $M^{4i+2}_\mathbb{C}$ contains a non zero class of type $(2i+1,2i+1)$, of the form
$\beta=\alpha\cup\overline{\alpha}$. On the other hand, $M^{2i+1,2i+1}_\mathbb{C}$ contains
$\omega^{2i+1}$. But these two cohomology classes
cannot be proportional because
$\beta^2=0$, while $\omega^{4i+2}\not=0$ by Lefschetz property and  because $4i+2\leq n$.
\cqfd
\subsection{Sub-Hodge structures and a lemma of Deligne}
The following lemma, communicated to the author by Deligne \cite{deligne}, and very much used in
\cite{voisin}, \cite{voisinjdg}, allows to detect sub-Hodge structures in a given cohomology algebra
endowed with a Hodge structure.
Let $A^*=\oplus_kA^k$ be a rational (resp. real) cohomology algebra endowed with a
Hodge structure. Let
$A^*_\mathbb{C}:=A^*\otimes\mathbb{C}$.
Let $Z\subset A^k_\mathbb{C}$ be an algebraic subset which is defined by
homogeneous equations expressed only using the ring structure on $A^*$. The examples we shall consider
in this paper  will often be of the form :

$$Z=\{\alpha\in A^k_\mathbb{C}/\alpha^l=0\,\,{\rm in }\,\,A^{kl}_\mathbb{C}\},$$
where $l$ is a given integer.

\begin{lemm}\label{dele} Let $Z$ be as above, and let $Z_1$ be an irreducible component of $Z$.
Assume the $\mathbb{C}$-vector space $<Z_1>$ generated by $Z_1$ is defined over $\mathbb{Q}$, (resp. over
$\mathbb{R}$),
that is $<Z_1>=B^k_\mathbb{Q}\otimes\mathbb{C}$ for some $B_\mathbb{Q}^k\subset A_\mathbb{Q}^k$
(resp. $<Z_1>=B^k_\mathbb{R}\otimes\mathbb{C}$ for some $B_\mathbb{R}^k\subset A_\mathbb{R}^k$).
Then $B_\mathbb{Q}^k$ (resp.  $B_\mathbb{R}^k$) is a rational
(resp. real) sub-Hodge structure of $A^k_\mathbb{Q}$ (resp. $A^k_\mathbb{R}$).
\end{lemm}
We refer to \cite{voisin} for the  proof of this lemma.
\section{Stability results \label{sec2}}
\subsection{Products}
We prove in this section   the following result:
\begin{theo} \label{product}Assume there is a polarized Hodge structure on a
connected (rational or real) cohomology algebra
$M$, and assume that
$$M\cong A\otimes B,$$
where $A$ and $B$ are (rational or real) cohomology algebras. If either
$A^1=0$ or $B^1=0$, then $A$ and $B$ are of even dimension
and there are polarized Hodge structures on
$A$ and $B$, inducing that of $M$.
\end{theo}
Geometrically, this implies that if a product $X\times Y$, where $X$ and $Y$ are smooth
compact oriented manifolds,  has the cohomology algebra  of a K\"ahler
compact manifold, and one of them has $b_1=0$, then
the cohomology algebras of  $X$ and $Y$ carry polarized Hodge structures, and thus
in particular
inherit all the constraints described in next section.

Another geometric consequence is the fact that if a K\"ahler compact  manifold $M$
is homeomorphic to a product $X\times Y$, where either $b_1(X)=0$ or $b_1(Y)=0$, then
any polarized Hodge structure on $H^*(M)$ comes from polarized Hodge structures
on $H^*(X)$ and $H^*(Y)$, a result which can be
compared to a result in deformation theory: if a
compact complex manifold is a product $X\times Y$, and $b_1(X)=0,\,b_1(Y)=0$, then small deformations
of $X\times Y$ are of the form
$X'\times Y'$, where $X'$ is a deformation of $X$ and
$Y'$ is a deformation of $Y$. (This is not true if only one of the assumptions
$b_1(X)=0,\,b_1(Y)=0$ is satisfied, as a product may deform to a non trivial fiber bundle.)

\begin{rema}{\rm The assumption that $A^1=0$ or $B^1=0$ is obviously necessary
in all the statements above, as the case of a complex torus shows. Indeed, the complex torus is a
product of an even number $2n$ of copies of $S^1$'s, and thus can be
written as a product of
$X$, a product of an odd number $2d-1$ of copies of $S^1$'s and
$Y$, a product of an odd number $2d'+1$ of copies of $S^1$'s, with $d+d'=n$.}
\end{rema}

To start the proof of the theorem, let us show the assertion of even dimensionality.
\begin{lemm}\label{ledim} Under the assumptions of Theorem \ref{product},
$A$ and $B$ are  even dimensional.
\end{lemm}
{\bf Proof.} Indeed, let $s=dim\,A$,
$t=dim\,B$, so that
$$s+t=2n:=dim\,M.$$
As $A^1=0$, or $B^1=0$, and $M=A\otimes B$,
one has
$$M^2=A^2\otimes B^0\oplus A^0\otimes B^2,$$
where we can identify canonically
$A^0,\,B^0$ with $\mathbb{Q}$.
Let $\omega\in M^2_\mathbb{R}$ be the class of a polarization. Then $\omega^n\not=0$ in
$M^{2n}=A^s\otimes B^t$.

Writing $\omega=a+b$, with $a\in A^2$ and $b\in B^2$, we conclude by writing
$$\omega^n=\oplus_{i\leq n} \binom{n}{i}a^i\otimes b^{n-i}$$  that
$A^s$ must be generated by a power of $a$, and $B^t$ must be generated by a power
of $b$. Thus $s$ and $t$ are even.

\cqfd
Note also that the proof showed that
the algebras $A$ and $B$ have their top degree part generated by a power of a degree
$2$ element.

Let us now apply Deligne's Lemma \ref{dele} to get the following:
\begin{lemm} \label{leh2} The subspaces
$$A^2\cong A^2\otimes B^0\subset M^2,\,B^2\cong A^0\otimes B^2\subset M^2$$
are rational sub-Hodge structures of $M^2$.
\end{lemm}
{\bf Proof.} By Deligne's lemma, it suffices to show how to recover these
subspaces algebraically, using only the algebra structure of $M$.
Let $2s=dim\,A,\,2t=dim\,B$, so that $s+t=n$.
We claim that
$A^2_ \mathbb{C}\subset M^2_\mathbb{C}$ is an irreducible component of
$$Z\subset M^2_\mathbb{C},\,Z=\{m\in M^2_\mathbb{C},\,m^{s+1}=0\},$$
and similarly for $B^2_\mathbb{C}$, with $s$ replaced by $t$.

To see this, let as before $\omega=a+b$ be a decomposition of a polarizing class $\omega\in M_\mathbb{R}$.
Then we proved in the previous lemma that
$a^s\not=0$ in $A^{2s}$.
The tangent space to $Z$ at the point $a$ is described as
$$T_{Z,a}=\{m\in M^2_\mathbb{C},\,a^sm=0\}.$$
Writing $m=\alpha+\beta,\,\alpha\in A^2,\,\beta\in B^2$, we conclude immediately
 that $\beta=0$ for $m\in T_{Z,a}$,
which shows that the Zariski tangent spaces of
$A^2_\mathbb{C} $ and $Z$ coincide at $a$.
As $A^2_\mathbb{C} \subset Z$  is smooth, this implies that $A^2_\mathbb{C}$ is an irreducible component
of $Z$.
\cqfd
Recall that a Hodge class in a rational Hodge structure $M^{2k}$ of
weight $2k$ is a rational element of
$M^{2k}$ which is also in $M^{k,k}$. As a corollary of the previous two lemmas, we get
\begin{coro} The $1$-dimensional rational spaces $ A^{2s}\subset M^{2s}$, $ B^{2t}\subset
M^{2t}$  are generated by Hodge
classes $\eta_A,\,\eta_B$ of $M$ of
respective  degrees
$2s,\,2t$.
\end{coro}
{\bf Proof.}
Indeed, consider the case of $A$. Then, as mentioned above, $A^{2s}$ is the image
of the map
$$Sym^sA^2\rightarrow Sym^sM^2\rightarrow M^{2s},$$ where
the first map is the inclusion, and the second is given by the product of $M$.
As $A^2$ is a sub-Hodge structure of $M^2$, $A^{2s}$ is also a sub-Hodge structure of
$M^{2s}$. As it is one dimensional, it must be generated by a Hodge class.
\cqfd
We finally prove the following:
\begin{lemm} \label{lefsch}Let $\omega\in M^2_\mathbb{R}$ be a polarizing class and decompose
it as
$$\omega=a+b,\,a\in A^2_\mathbb{R},\,b\in B^2_\mathbb{R}.$$ Then the classes
$a$ and $b$ satisfy the Lefschetz property, namely
for any $k\leq s$
$$\cup a^{s-k}:A^k_\mathbb{R}\rightarrow A^{2s-k}_\mathbb{R}$$
is an isomorphism, and similarly for $b$ and $B$.
\end{lemm}
{\bf Proof.} Let $\alpha\in A^k_\mathbb{R}\subset M^k_\mathbb{R}$ and assume that
$a^{s-k}\cup\alpha=0$ in $A^{2s-k}_\mathbb{R}$. It then follows that
$a^i\cup\alpha=$ in $A^{k+2i}_\mathbb{R}$ for all $i\geq s-k$.
Let us compute now:
$$\omega^{n-k}\cup\alpha=\sum_{i\leq n-k}\binom{n-k}{i} a^i\cup b^{n-k-i}\cup\alpha.$$
As $a^i\cup\alpha=0$ for $i\geq s-k$, the sum runs only over the pairs
$(i,j)$ with $i+j=n-k,\,i<s-k$, and thus $j>n-s=t$. As $b^{t+1}=0$, we get
$\omega^{n-k}\cup\alpha=0$ and the Lefschetz property for $(M,\omega)$ shows that
$\alpha=0$.
\cqfd
{\bf Proof of theorem \ref{product}.} We assert that each
$A^i,\,B^i\subset M^i$ is a sub-Hodge structure of $M^i$.
For $i=1$, this follows from the fact that we have either
$A^1=0,\,B^1=M^1$ or $A^1=M^1,\,B^1=0$ and for $i=2$, this follows from
Lemma \ref{leh2}. We choose now a polarizing class
$\omega\in M^{1,1}_\mathbb{R}$ and decompose it into
$\omega=a+b$. Then by Lemma \ref{leh2}, $a$ and $b$ are of type $(1,1)$.

We prove now the result by induction on $i$. So assume the result is proved
for $i-1$. Then all the
subspaces
$$A^l\otimes B^{l'},\,l+l'=i,\,l>0,l'>0$$ of $M^{i}$
are sub-Hodge structures. Using now Lemma \ref{leh2}, and using a polarization
$\omega=a+b\in M^{1,1}_\mathbb{R}$,  we have that
$a$ and $b$ are in $M^{1,1}_\mathbb{R}$ and thus we get that
each subspace
$$a^{s-l}\cup b^{t-l'}\cup A^l_\mathbb{R}\otimes B^{l'}_\mathbb{R}\subset M^{2n-i}_\mathbb{R}$$
is a real sub-Hodge structure.

By lemma \ref{lefsch}, this subspace is in fact equal to
$A^{2s-l}\otimes B^{2t-l'}\otimes \mathbb{R}$, and thus we conclude that
for $l>0$ and $l'>0$, $l+l'=i$
$$A^{2s-l}\otimes B^{2t-l'}\subset M^{2n-i}$$ is a rational sub-Hodge structure
of $M^{2n-i}$, because it is rational, and tensored
by $\mathbb{R}$ becomes a real sub-Hodge structure. But the orthogonal
of $$\oplus_{l>0,l'>0,l+l'=i}A^{2s-l}\otimes B^{2t-l'}$$ with respect
to the intersection pairing on $M$ is equal to
$A^i\oplus B^i$.
Thus $A^i\oplus B^i$ is a sub-Hodge structure of $M^i$.
But $A^i\subset A^i\oplus B^i$ is the kernel of the restriction
to $A^i\oplus B^i\subset M^i$ of
the multiplication by the Hodge class
$\eta_A: M^i\rightarrow M^{i+2s}$, where
$\eta_A$  generates $A^{2s}$, and similarly
$B^i\subset A^i\oplus B^i$ is the kernel of the restriction
to $A^i\oplus B^i\subset M^i$ of
the multiplication by the Hodge class
$\eta_B: M^i\rightarrow M^{i+2t}$, where
$\eta_B$  generates $B^{2t}$. Thus we conclude that $A^i$ and $B^i$ are sub-Hodge structures
of $M^i$.

To conclude the proof of the theorem, it remains to show that
$a$ polarizes the Hodge structure on $A\subset M$, and $b$ polarizes the Hodge structure
on $B\subset M$. We already proved that $a$ and $b$ satisfy the Lefschetz property.
We have to show that the second bilinear relations (\ref{secondriemann}) hold.

This is easy because, the same argument as in the proof of Lemma \ref{lefsch} shows that if
$\alpha\in A^i_\mathbb{R},\,i\leq s$ is primitive for $a$, namely satisfies
$a^{s-i+1}\cup\alpha=0$, then
$\alpha\in M^i_\mathbb{R}$ is primitive for $\omega$. Furthermore, we have, for primitive
$\alpha,\,\beta$, and for the choices of generators
$$\omega^n,\,a^s,\,b^t$$
of $M^{2n}_\mathbb{R},\,A^{2s}_\mathbb{R},\,B^{2t}_\mathbb{R}$ respectively:
$$\nu<\alpha,a^{s-i}\cup\beta>_A=<\alpha,\omega^{n-i}\cup\beta>_M,$$
where
$$\nu=\frac{\binom{n-i}{s-i}}{\binom{n}{s}}.$$
Thus the second bilinear relations for $(M,\omega)$ imply the second bilinear relations
for $(A,a)$.
\cqfd
Note to conclude that, without polarizations,
 the same arguments prove the following result, which will be used
later on:
\begin{theo}\label{product2} Assume there is a  Hodge structure on a
 (rational or real) cohomology algebra
$M$, and assume that
$$M\cong A\otimes B,$$
where $A$ and $B$ are (rational or real) cohomology algebras. If
$A^1=0$, and $A$ and $B$ are generated in degree $\leq 2$,
then $A$ and $B$ are of even dimension
and there are  Hodge structures on
$A$ and $B$, inducing that of $M$.
\end{theo}

{\bf Proof.} As in the previous proof, we first use Deligne's Lemma
\ref{dele} to prove that $A^2\otimes 1_B$ and $1_A\otimes B^2$ are sub-Hodge structures
of $M^2$. To prove this, note that $M$ is generated in degree $\leq 2$, which implies, as it is of even dimension
$2n$,
that there is a $\omega\in M^2$ such that $\omega^{n}\not=0$ in $M^{2n}$ (otherwise any monomial
$\omega_1\ldots\omega_n$ with $deg\,\omega_i=2$ would vanish in $M^{2n}$, contradicting the fact that
$M^{2n}$ is generated by a product of elements of degree $1$ or $2$).
Having this, and writing $\omega=\alpha+\beta,\,\alpha\in A^2\otimes 1_B,\,\beta\in 1_A\otimes B^2$,
we conclude as before that
$$\alpha^d\not=0\,\,{\rm in}\,\,A^2,\,\beta^{d'}\not=0\,\,{\rm in}\,\,B^2,$$
for some integers $d,\,d'$ such that $dim\,A=2d$, $dim\,B=2d'$.
It follows then that $A^2_{\mathbb{C}}\otimes 1_B$ is recovered as an irreducible component of
the set
$$\{m\in M^2_{\mathbb{C}},\,m^{d+1}=0\,\,{\rm in}\,\,M^{2d+2}_\mathbb{C}\},$$
and similarly for $1_A\otimes B^2_\mathbb{C}$. Thus by Deligne's Lemma \ref{dele},
$A^2\otimes 1_B$ and $1_A\otimes B^2$ are sub-Hodge structures
of $M^2$.

As $A^1=M^1$ the same is true in degree $1$. As the algebras are generated in degree
$\leq2$, it follows that $A=A\otimes 1_B$ and $B=1_A\otimes B$ are sub-Hodge structures of $M$.
\cqfd
\subsection{Projective bundles}

We consider now the simplest kinds of compact oriented manifolds which are close to be a product:
namely
bundles over a basis satisfying the Leray-Hirsch condition. Then
 their cohomology is additively the tensor product of the cohomology of the basis
and the cohomology of the fiber, but not multiplicatively. The simplest examples of this are given
by complex projective bundles $\mathbb{P}(E)\rightarrow X$ associated to complex vector bundles on
compact oriented manifolds.

\begin{theo} \label{theoprojbun} Let $X$ be a compact smooth oriented manifold,  and let
$E$ be a complex vector bundle of rank $\geq2$ on $X$ with trivial determinant. Assume that the cohomology
of $X$ is generated in degrees $1$ and $2$. Then if the cohomology algebra
$H^*(\mathbb{P}(E),\mathbb{Q})$ carries a  Hodge structure, the cohomology algebra
$H^*(X,\mathbb{Q})\subset H^*(\mathbb{P}(E),\mathbb{Q})$ has an induced Hodge  structure, for which
the Chern classes $c_i(E)$ are Hodge classes. A similar result holds with
$H^*(\mathbb{P}(E),\mathbb{Q}), H^*(X,\mathbb{Q})$
replaced by $H^*(\mathbb{P}(E),\mathbb{R}),\,H^*(X,\mathbb{R})$.
\end{theo}
{\bf Proof.} Let $X$ be of real dimension $2n$,
 $\pi:\mathbb{P}(E)\rightarrow X$ be the structural map. We know that
$$\Gamma:=\pi^*H^2(X,\mathbb{Q})\subset H^2(\mathbb{P}(E),\mathbb{Q})$$ is a hyperplane, which satisfies
the properties that
$\alpha^{n+1}=0$, for any $\alpha\in \Gamma$.
Next, the cohomology of $\mathbb{P}(E)$, as the cohomology of $X$, being generated in degree
$1$ and $2$, it follows that there exists a class
$\beta\in H^2(\mathbb{P}(E),\mathbb{Q})$ such that $\beta^{n+r-1}\not=0$, where $r:=rank\,E$.

Thus we conclude that the hyperplane
$\Gamma_{\mathbb{C}}\subset  H^2(\mathbb{P}(E),\mathbb{C})$ must be
an irreducible component of the algebraic subset
$$Z:=\{\beta\in H^2(\mathbb{P}(E),\mathbb{C}),\,\beta^{n+r-1}=0\}.$$
By Deligne's Lemma \ref{dele}, it follows that $\Gamma$ is a sub-Hodge structure of
$H^2(\mathbb{P}(E),\mathbb{Q})$ for the given Hodge structure.
As $\pi^*:H^1(X,\mathbb{Q})\rightarrow H^1(\mathbb{P}(E),\mathbb{Q})$ is an isomorphism, we
have the same conclusion for the cohomology of  degree $1$. Finally, as the cohomology of
$X$ is generated in degrees $1$ and $2$, we conclude that the cohomology sub-algebra
$$\pi^*H^*(X,\mathbb{Q})\subset H^*(\mathbb{P}(E),\mathbb{Q})$$
is also a sub-Hodge structure.

Observe now that the injective map
$$\pi^*:H^*(X,\mathbb{Q})\rightarrow H^*(\mathbb{P}(E),\mathbb{Q})$$
admits as its dual map the Gysin map
$$\pi_*:H^*(\mathbb{P}(E),\mathbb{Q})\rightarrow H^{*-2r+2}(X,\mathbb{Q})$$
which thus must be also a morphism of Hodge structures (of bidegree $(-r+1,-r+1)$ on each graded piece).

We claim that there is a  class $\beta\in H^{2}(\mathbb{P}(E),\mathbb{Q})$, unique up
to a multiplicative coefficient, such that $\beta$ does not vanish modulo
$\pi^*H^2(X,\mathbb{Q})$ and satisfies
$$\pi_*\beta^{r}=0\,\,{\rm in}\,\,H^2(X,\mathbb{Q}).$$
Furthermore this class must be a Hodge class.
Indeed, we take for $\beta$ the class $c_1(H)$ where $H$ is the dual of the
Hopf line bundle on $\mathbb{P}(E)$. Then it is a standard fact
(see \cite{fulton}) that
$$\pi_*\beta^{r}=- c_1(E).$$
As $c_1(E)$ was supposed to be $0$, this proves the existence of $\beta$. As for the uniqueness, observe
 that for $\alpha\in H^2(X,\mathbb{Q})$,
$$\pi_*(\beta+\pi^*\alpha)^r=\pi_*\beta^r+r\alpha,$$
so that $\pi_*(\beta+\pi^*\alpha)^r=0$ implies $\alpha=0$.

The same argument shows that the complex line $\beta\mathbb{C}\subset H^2(\mathbb{P}(E),\mathbb{C})$
is an irreducible component of the closed algebraic subset
$$Z:=\{\gamma\in H^2(\mathbb{P}(E),\mathbb{C}),\,\pi_*\gamma^{r}=0\}.$$
As $\pi_*$ is a morphism of Hodge structures, we conclude by Deligne's Lemma \ref{dele}
that $\beta$ must be a Hodge class.

The proof is now finished. Indeed, let $\beta$ be defined (up to a multiplicative coefficient)
 as above. Then
$\beta$ satisfies  in $H^*(\mathbb{P}(E),\mathbb{Q})$ a unique polynomial equation
$$\beta^{r}=\sum_{1\leq i\leq r} \beta^{r-i}\cup\pi^*\alpha_i,$$
where the $\alpha_i$ are proportional to $c_i(E)$. As all the powers  $\beta^l$'s are Hodge classes
on $\mathbb{P}(E)$, and the $\alpha_i$'s can be recovered
as polynomials in the Segre classes
$\sigma_j(E)=\pi_*(\beta^{r-1+j})$ (cf \cite{fulton}),
 the $\alpha_i$ must be also Hodge classes on $X$ for the induced
Hodge structure on $H^*(X,\mathbb{Q})$.

\cqfd
\section{Explicit constraints\label{sec3}}
The purpose of this section is to describe explicit constraints on a cohomology algebra, imposed by
the presence of a (polarized) Hodge structure. We want to separate the constraints
on the {\it real} cohomology algebra imposed by the Hodge decomposition, which will be
considered in section \ref{realsec}, from more subtle
constraints related to the {\it rational} cohomology algebra, which will be explained in
section \ref{qsec} and from those imposed by the polarisation (section \ref{sectionpol}).
We will illustrate  the effectiveness of each of these criteria by exhibiting compact symplectic manifolds
satisfying the Lefschetz property and not satisfying the considered criterion.
\subsection{Constraints coming from the real Hodge structure
\label{realsec}}
It is clear that constraints coming only from the real (or rational) Hodge structure, without polarizations,
must involve the odd degree cohomology. Indeed, if there is no odd degree cohomology, we can
put the trivial Hodge structure on all the cohomology groups, and this will  obviously satisfy
the compatibility conditions (\ref{compatcond}).

The classically known topological restriction on the cohomology algebra of a compact K\"ahler
manifold coming from the Hodge decomposition is the fact that
odd Betti numbers $b_{2i+1}$ must be even.
This condition is very restrictive for surfaces, as it is known that
among compact complex surfaces, this condition characterizes
the K\"ahler ones, a result due to Kodaira \cite{kodaira}.

Let us refine this restriction using the cohomology algebra.
\begin{lemm}\label{basic} Let $M$ be a real cohomology algebra. If $M$ carries a  Hodge structure,
then for any pair of integers $k,\,l$ with $k$ even and $l$ odd, the  product:
$$\mu: M^l\otimes M^k\rightarrow M^{k+l}$$
must be of even rank. More generally, if $M^{'k}\subset M^k,\,M^{'l}\subset M^l$
 are sub-Hodge structures, then
the  product:
$$\mu: M^{'l}\otimes M^{'k}\rightarrow M^{k+l}$$
must be of even rank.
\end{lemm}
{\bf Proof.} Indeed, as $M^{'k}\otimes\mathbb{C}$ and $M^{'l}\otimes\mathbb{C}$
are stable under the Hodge decomposition, that is decompose into the
direct sum of their components of type $(p,q),\,p+q=k$, resp. $p+q=l$, it follows
from  the compatibility conditions
(\ref{compatcond}) that the image $Im\,\mu_\mathbb{C}$ of $\mu$, tensored by $\mathbb{C}$, is
also stable under the Hodge decomposition,
namely, it is the direct sum of its terms of type $(p,q),\,p+q=k+l$. As
$k+l$ is odd, and $Im\,\mu_\mathbb{C}$ is defined over
$\mathbb{R}$, the Hodge symmetry (\ref{hodgesymmetry}) is satisfied
by $Im\,\mu_\mathbb{C}$, which implies that it is of even rank.
\cqfd
This lemma can be combined with Deligne's Lemma \ref{dele}, to get effective
restrictions which are much stronger than the classical ones.
Let us state this restriction   explicitely:
\begin{prop} \label{utilepourcequisuit}Let $M$ be a real cohomology algebra with Hodge structure. For some even
integer $2k$, let $Z_1,\ldots, Z_r\subset M^{2k}_\mathbb{C}$ be  algebraic subsets
as in Lemma \ref{dele}. Suppose the complex vector spaces $<Z_i>,\,i=1,\ldots,\,r$
 are defined over $\mathbb{R}$, that is
$$<Z_i>=B_i\otimes\mathbb{C},\,B\subset M^{2k}.$$
Let $B:=\sum_iB_i\subset M^{2k}$.
Then for any odd integer $l$, the product map
$$ B\otimes M^l\rightarrow M^{2k+l}$$
has even rank.

\end{prop}
Indeed, Deligne's Lemma then tells us that $B$ has to be a sub-Hodge structure, so that we can apply Lemma
\ref{basic}.\cqfd
 Let us construct using
Theorem \ref{utilepourcequisuit} a compact symplectic manifold $X$ which satisfies the Lefschetz
condition, hence in particular has its odd Betti numbers even, but whose {\it real}
cohomology algebra
does not admit any Hodge structure. So in particular, this
$X$ does not have the real cohomology algebra of a compact K\"ahler manifold.

\begin{ex} We start with a a complex torus $T=\mathbb{C}^3/\Gamma$ of dimension $3$, and we fix a symplectic structure
on $T$, given for example by a constant K\"ahler form $\omega$
on $\mathbb{C}^3$. We may even assume that $T$ is an abelian variety
and that the cohomology class
$\lambda_0$ of the K\"ahler form is rational.
Next,  we choose
$2$ elements $\lambda_1,\,\lambda_2$ in $H^2(T,\mathbb{Q})$ which satisfy the property
that
$$\lambda_1\cup H^1(T,\mathbb{Q})+\lambda_2\cup H^1(T,\mathbb{Q})$$
has rank $11$. As $rank\,H^1(T,\mathbb{Q})=6$, this means that
the map
\begin{eqnarray}\label{mu}\mu:H^1(T,\mathbb{Q})\oplus H^1(T,\mathbb{Q})\rightarrow H^3(T,\mathbb{Q}),\\ \nonumber
\mu((a,b))=\lambda_1\cup a+\lambda_2\cup b
\end{eqnarray}
has a $1$-dimensional kernel.
Let us give an explicit example of such a pair:
Choose a basis $w_1,\ldots,w_6$ of the $\mathbb{Q}$-vector space
$H^1(T,\mathbb{Q})\cong \Gamma^*_\mathbb{Q}$.
The cohomology $H^*(T,\mathbb{Q})$ identifies to the exterior algebra $\bigwedge\Gamma^*_\mathbb{Q}$,
the cup-product being identified with the exterior product.
Let
$$\lambda_1=w_1\wedge w_2+w_3\wedge w_4,\,\lambda_2=w_3\wedge w_5-w_1\wedge w_4.$$
Then we clearly have $w_1\wedge \lambda_1=w_3\wedge\lambda_2$, and it is an easy exercise
to verify that the kernel of the map (\ref{mu}) is generated by this relation.

We now do the following : choosing
$\epsilon \in\mathbb{Q}$ small enough, the classes $\lambda_0+\epsilon \lambda_1$
and $\lambda_0+\epsilon \lambda_2$ are symplectic classes, which can be represented by symplectic forms
$\omega_1,\,\omega_2$
in the same deformation class as $\omega$. In fact, the important point for us is
the fact that $\omega_0+\omega_1+\omega_2$ is a symplectic form, whose class is close to
$3\lambda_0$.

As the classes $\lambda_0,\,\lambda_1,\,\lambda_2$ are rational, we can find
multiples $M\lambda_0,\,M\lambda_1,\,M\lambda_2$ which consist of integral symplectic classes,
such that there are symplectic embeddings (we can use
\cite{gromov}, 3.4, or  approximately holomorphic  embeddings (see \cite{Munoz} using results of \cite{donaldson}))
$$\phi_i: T\rightarrow \mathbb{P}^N,\,i=1,\,2,\,3,$$
with $\phi_i^*\Omega=M\omega_i$, where $\Omega$ is the Fubini-Study symplectic form on $\mathbb{P}^N$.

Let $\psi:T\rightarrow \mathbb{P}^N\times\mathbb{P}^N\times \mathbb{P}^N$ be the
map $(\phi_1,\,\phi_2,\,\phi_3)$.
For the product symplectic form
$\widetilde{\Omega}=p_1^*\Omega+p_2^*\Omega+p_3^*\Omega$, the image $W=\psi(T)$ is a symplectic submanifold,
because $\psi^*\widetilde{\Omega}=M(\omega_0+\omega_1+\omega_2)$.

Our example will be the symplectic manifold $X$ obtained as the symplectic blow-up of
$\mathbb{P}^N\times\mathbb{P}^N\times \mathbb{P}^N$ along $W$.
For an adequate choice of symplectic form
of class
$\mu_1+\mu_2+\mu_3-\eta e$, with $\eta$ very small, the Lefschetz property is satisfied,
as it follows easily from the fact that
the restriction of the  symplectic class $p_1^*\omega+p_2^*\omega+p_3^*\omega$ to
$W=T$ satisfies the Lefschetz property on $T$.
(The Lefschetz property for symplectic blow-ups is studied in general in
\cite{caval}.)
\end{ex}

\begin{prop}\label{rangpair} The cohomology algebra of $X$
 does not satisfy
the condition of
Proposition \ref{utilepourcequisuit}, hence does not admit any real Hodge structure. In particular
$X$ does not have the cohomology algebra of a compact K\"ahler manifold.
\end{prop}
{\bf Proof.}
By the computation of the cohomology of a symplectic blow-up (cf \cite{voisinbook}, 7.3.3
in the complex case), we get
that
$H^2(X,\mathbb{Q})=\mathbb{Q}^4$, generated by the class $e$ of the exceptional divisor,
and the classes $\mu_i:=\tau^*(p_i^*[\Omega]),\,i=1,\,2,\,3$, where
$\tau:X\rightarrow \mathbb{P}^N\times\mathbb{P}^N\times \mathbb{P}^N$
is the blowing-up map. Furthermore, letting
$j:E\rightarrow X$ denote the inclusion of the exceptional divisor
of $\tau$, and $\tau':E\rightarrow W=T$ the restriction of $\tau$, we have that
$H^3(X,\mathbb{Q})=j_*(\tau^{'*}H^1(T,\mathbb{Q}))$.

Now, one sees easily that each class $\mu_i$ generates an irreducible component of the algebraic subset
$$Z=\{a\in H^2(X,\mathbb{C}),\,a^{N+1}=0\,\,{\rm in}\,\,H^{2N+2}(X,\mathbb{C})\}.$$
As each $\mu_i$ is rational, it must be a Hodge class by Deligne's Lemma \ref{dele}.

 Hence we proved
that for any Hodge structure on $H^*(X,\mathbb{Q})$,  the classes
$\mu_i$ are of type $(1,1)$. It remains to see
that for two adequately chosen rational (or real) cohomology classes
$\lambda,\lambda'$ which are combinations of the $\mu_i$,
 the rank of
\begin{eqnarray}\label{autremu}\mu:H^3(X,\mathbb{Q})\oplus H^3(X,\mathbb{Q})\rightarrow H^5(X,\mathbb{Q})\\
\nonumber (\alpha,\beta)\mapsto \lambda\cup\alpha+\lambda'\cup\beta
\end{eqnarray}
is odd.

Let us take
$$\lambda=\mu_2-\frac{1}{\epsilon}\mu_1,\,\lambda'=\mu_3-\frac{1}{\epsilon}\mu_1.$$
Then as $\lambda$, resp. $\lambda'$, is the pull-back via $\tau $ of a class
$\tilde{\lambda}$ on
$\mathbb{P}^N\times\mathbb{P}^N\times \mathbb{P}^N$,
it follows that the following diagram commutes
$$\begin{matrix} &\tilde{\lambda}_{\mid W}\cup:&H^1(W,\mathbb{Q})&\rightarrow& H^3(W,\mathbb{Q})\\
&&j_*\tau^{'*}\downarrow&& j_*\tau^{'*}\downarrow\\
&{\lambda}\cup:&H^3(X,\mathbb{Q})&\rightarrow &H^5(X,\mathbb{Q}),
\end{matrix}
$$
where the vertical maps are injective. We have the similar result for $\lambda'$.
On the other hand,  identifying $W$ with $T$, we have
$$\tilde{\lambda}_{\mid W}=M\lambda_1,\,\,\tilde{\lambda'}_{\mid W}=M\lambda_2.$$

As the map $j_*\tau^{'*}:H^3(T,\mathbb{Q})\rightarrow H^5(T,\mathbb{Q})$ is
injective, the fact that the map
(\ref{mu}) has odd rank  thus implies that the map (\ref{autremu}) has odd rank.

\cqfd

Let us now describe another non trivial necessary condition for a real cohomology algebra
$M$
to admit a Hodge structure.
Let $H^{2i+1}$ be an even dimensional  real  vector space, and
let $H^{4i+2}$ be a real vector space. Let
$$\mu: {\bigwedge}^2H^{2i+1}\rightarrow H^{4i+2}$$
be a linear map.
\begin{lemm}\label{adhoc}  Suppose there are Hodge
structures of respective weights $2i+1$ and $4i+2$ on
$H^{2i+1}$ and $H^{4i+2}$ for which $\mu$ is a morphism of Hodge structures. Then there exists a complex vector subspace
$W\subset H^{2i+1}_\mathbb{C}$ satisfying the properties:
\begin{enumerate}
\item \label{item1juin}
$rk\,W=\frac{1}{2} rk\,H^{2i+1}_\mathbb{C}$.
\item $rk\,\mu_{\mid \bigwedge^2W}\leq \frac{1}{2} rk\,\mu$.
\end{enumerate}
\end{lemm}
{\bf Proof.} Indeed,
 let
$$W:=F^iH^{2i+1}_\mathbb{C}:=H_\mathbb{C}^{2i+1,0}\oplus\ldots\oplus H_\mathbb{C}^{i+1,i}.$$
Then by Hodge symmetry, \ref{item1juin}
is satisfied. Furthermore, if $\mu$ is a morphism of
Hodge structures, we get $$\mu({\bigwedge}^2F^iH^{2i+1}_\mathbb{C})\subset
H_\mathbb{C}^{4i+2,0}\oplus\ldots\oplus H_\mathbb{C}^{2i+2,2i},$$
and as $Im\,\mu$ is a real sub-Hodge structure of $H^{4i+2}$, we conclude
by Hodge symmetry for $Im\,\mu$ that we have the condition
\begin{eqnarray}\label{rang12} rank\,\mu({\bigwedge}^2F^iH^{2i+1}_\mathbb{C})\leq \frac{1}{2}rank\,\mu.
\end{eqnarray}
\cqfd
Coming back to real cohomology algebras, we can apply this lemma to
$H^{2i+1}=M^{2i+1}$, where $\mu$ is  the
cup-product map. Moreover, combining  this
lemma  with the results of the previous section,
we can also apply  it
to more general morphisms of Hodge structures of the form:
$$\mu'=\phi\circ \mu:{\bigwedge}^{2}M^{2i+1}\rightarrow M^{4i+2}\stackrel{\phi}{\rightarrow}
M',$$
where $M'$ will be an adequate Hodge structure of even
weight and $\phi$ will be shown to be a morphism of Hodge structures.

In order to show the effectiveness of this restriction on the structure of
cohomology algebras,
we need the following lemma:
\begin{lemm}\label{genericgrass} Let $M$ be a  real or rational vector space
of rank $2n$,
and let $$\mu:{\bigwedge} ^2M\rightarrow Q$$ be a generic linear  surjective map
to a rational or real vector space $Q$ of rank  $q=2q'$ or $q=2q'+1$ satisfying
the following numerical conditions :
\begin{eqnarray}\label{condnum} q'\leq \frac{n(n-1)}{2},\,(q-q')(\frac{n(n-1)}{2}-q')>n^2.
\end{eqnarray}
Then there is no complex subspace
$W\subset M_\mathbb{C}$
of rank $n$, such that
\begin{eqnarray}\label{*}rank\,\mu_{\mid\bigwedge^2W}\leq q'.\end{eqnarray}
\end{lemm}
{\bf Proof.} We work inside the space
$K=Hom\,(\bigwedge ^2M_\mathbb{C},Q_\mathbb{C})$ of $\mathbb{C}$-linear maps
$\mu:\bigwedge ^2M_\mathbb{C}\rightarrow Q_\mathbb{C}$ and make a dimension count.
The dimension of the Grassmannian $Grass(n,M_\mathbb{C})$ of $n$-dimensional subspaces
of $M_\mathbb{C}$ is $n^2$. For fixed $W\in Grass(n,M_\mathbb{C})$, the codimension of
the (Zariski closed) space consisting of $\mu\in K$
satisfying the condition that $rank\,\mu_{\mid\bigwedge^2W}\leq q'$
is equal to $(q-q')(\frac{n(n-1)}{2}-q')$ assuming this number is $\geq 0$,
which is implied by our first assumption. It follows immediately
that under our first assumption in (\ref{condnum}), the codimension of the set of
$\mu$ for which there exists a $W$ such that (\ref{*}) holds
is at least $q'(\frac{n(n-1)}{2}-q')-n^2$. This is positive under our second assumption
in (\ref{condnum}),
and it follows that this  Zariski closed subset of $K$ is a proper subset.
\cqfd

Let us now construct an example of a symplectic compact manifold
$X$  whose real cohomology algebra does not satisfy
the  criterion given by Lemma \ref{adhoc}.

\begin{ex} In lemma \ref{genericgrass}, we make $n=5$, $q'=5$, $q=2q'+1$.
We now consider a real torus $T$ of dimension $10$, and
let $K$ be any simply connected compact K\"ahler manifold
satisfying the condition that $rank\, H^2(K,\mathbb{Q})= 11$ and that the cohomology
of $K$ is generated in degree $2$.

We consider now $T\times K$, which can be endowed
with the structure of a compact K\"ahler manifold, hence in particular
is  a compact symplectic
 manifold. We have,
using K\"unneth decomposition and Poincar{\'e} duality, an inclusion

\begin{eqnarray}\label{kunneth}
Hom\,(H^2(T,\mathbb{Q}),H^2(K,\mathbb{Q}))\subset H^{10}(T\times K,\mathbb{Q}).
\end{eqnarray}
Choose now a generic surjective map
$$\mu:{\bigwedge}^2H^1(T,\mathbb{Q})=H^2(T,\mathbb{Q})\rightarrow H^2(K,\mathbb{Q})$$
and let $\lambda\in H^{10}(T\times K,\mathbb{Q})$ be the class which is the image of $\mu$
under the map (\ref{kunneth}).

Our example will be the symplectic manifold
$\mathbb{P}(E)$, for any complex vector bundle $E$ on $Y:=T\times K$
such that
$$c_1(E)=0, \,c_5(E)=m\lambda$$
for some non zero integer $m$.
\end{ex}
Let us now  combine Theorem \ref{theoprojbun} and Lemma \ref{adhoc} to show the
following result.
\begin{prop} \label{theo11mars}
 The compact symplectic manifold
 $\mathbb{P}(E)$ has  the property that there is no Hodge structure on the
real cohomology algebra
$H^*(\mathbb{P}(E),\mathbb{R})$.

\end{prop}

{\bf Proof.} Suppose the conclusion
is not satisfied. We first use Theorem \ref{theoprojbun} (version with real coefficients) to conclude that
under our assumptions, there should be a real Hodge structure on
$H^*(Y,\mathbb{R})=H^*(T\times K,\mathbb{R})$ for which $\lambda$ is a Hodge class.

We next use Theorem \ref{product2}, (version with real coefficients) to conclude that
there are real Hodge structures on $H^*(T,\mathbb{R})$ and $H^*(K,\mathbb{R})$ which induce the
 Hodge structure on
$H^*(T\times K,\mathbb{R})$.
The class $\lambda$ being a Hodge class on $T\times K$, we conclude that the corresponding morphism
$$\mu:{\bigwedge}^2H^1(T,\mathbb{R})\rightarrow H^2(K,\mathbb{R})$$ is a morphism
of Hodge structures.
But this contradicts the fact that $\mu$ is generic, so that by Lemma
\ref{genericgrass}, the map $\mu$ does not satisfy the necessary condition
(\ref{rang12}), with $2i+1=1$ in this case.

\cqfd

\begin{rema}{\rm The  two criteria we applied in this section
 to detect the non-existence
of real Hodge structures on the cohomology algebras of
certain symplectic manifolds, and the proof of their
effectiveness, are not only  different geometrically.
 The difference of nature is made clear in the fact
that in the previous case, our criterion did not apply to a
generic choice of $\lambda_1,\,\lambda_2$, and indeed, for a generic choice of
$\lambda_0, \,\lambda_1,\,\lambda_2$, there exists a rational Hodge structure on the cohomology
of the variety we constructed.

In the second example, a {\it generic} choice of $\mu$ will lead to an example
 where the criterion applies.}
\end{rema}
\subsection{Constraints coming from the rational Hodge structure \label{qsec}}
In this section, we want to use Theorem \ref{theoprojbun} to construct a
compact symplectic manifold $X$, whose
cohomology algebra does not carry any Hodge structure, polarized or not.
Working a little more, one can
see
that the {\it rational} information is actually needed, the information given by the real cohomology
algebra being too weak.

Here, it is hard to formulate the criterion in a general way. The starting point is however
Theorem \ref{theoprojbun}, which tells us that for some types of manifolds,
Hodge structures on their cohomology algebras come from Hodge structures on
a certain
subalgebra, for which certain classes must be Hodge classes. This point can be combined with the
general observation made in \cite{voisin} that a rich cohomology algebra may prevent a cohomology
algebra with Hodge structure to carry any non trivial Hodge class.

We will content ourselves illustrating the combination of the two arguments on an example.

\begin{ex} We start with  the simplest example $Y$
  constructed in \cite{voisin} of compact K\"ahler manifold
  not having the rational cohomology algebra of a projective complex manifold, namely,
  we consider a complex torus
$$T={\Gamma}_\mathbb{C}/(\Gamma^{1,0}\oplus \Gamma),$$
where $\Gamma$ is a lattice of even rank $2n$; we assume there is an endomorphism   $\phi$
 acting on $\Gamma$, and that $\Gamma^{1,0}\subset \Gamma_\mathbb{C}$ is the eigenspace
 associated to the choice of $n$ complex eigenvalues of $\phi$, not pairwise
 conjugate. Here $\phi$ is assumed to have only complex eigenvalues, and to satisfy
 the following condition :

 \begin{eqnarray}
 \label{condition*}{\rm The\,\, Galois\,\, group\,\, of
 \,\, the \,\,splitting\,\, field\,\, of \,\,}\mathbb{Q}[\phi]\,\,
 {\rm is\,\, the }\\ \nonumber{\rm symmetric\,\, group\,\, in} \,\,2n \,\,{\rm letters\,\, acting
 \,\, on\,\, the }\,\,2n
 \,\,{\rm eigenvalues\,\, of\,\,
 }\phi.
 \end{eqnarray}

 The torus $T$ admits then the endomorphism $\phi_T$ induced by the $\mathbb{C}$-linear extension
 of $\phi$ acting on ${\Gamma}_\mathbb{C}$, which preserves
 $\Gamma^{1,0}$ and $\Gamma$.

 The K{\"a}hler manifold  $Y$ was obtained by successive blow-ups of $T\times T$.
 Namely, observing that the four subtori
 $$T_1=T\times 0,\,T_2=0\times T,\,T_3=Diag(T),\,T_4=Graph(\phi_T)$$
 of $T\times T$ meet pairwise transversally in finitely many points $x_1,\ldots , x_N$, we first blow-up
 the $x_i$'s, getting
 $\widetilde{T\times T}_{x_1,\ldots,x_N}$; then the proper transforms $\widetilde{T}_i$
 are smooth and disjoint,
 and we get $Y$ by blowing them in $\widetilde{T\times T}_{x_1,\ldots,x_N}$.

We shall denote $\tau:Y\rightarrow T\times T$ the natural map, which is the composition of two blow-ups.

 Let us choose now any integral cohomology class
 $\lambda\in H^2(Y,\mathbb{Z})$. Then $\lambda=c_1(L)$, where $L$ is a $C^\infty$ complex line bundle
 on $Y$. Our example will be  the manifold
 $$X=\mathbb{P}(L\oplus L^{-1}).$$

  \end{ex}
Let us show:
\begin{theo}\label{theotorus} If $dim\,T\geq 3$ and
$0\not=\lambda^2$ belongs to $\tau^*H^4(T\times T,\mathbb{Q})\subset H^*(Y,\mathbb{Q})$, then the
rational cohomology algebra of the manifold $X$ does not admit a Hodge structure.
\end{theo}
{\bf Proof.} We first observe  that the manifold $Y$ has its cohomology generated in
degrees $1$ and $2$. Indeed, this is true for the torus $T\times T$, (for which degree $1$ suffices),
thus also for $\widetilde{T\times T}_{x_1,\ldots,x_N}$ by
adding the classes of exceptional divisors over points. To see that this remains true for $Y$, observe
that the restriction maps in cohomology
$$H^*(\widetilde{T\times T}_{x_1,\ldots,x_N},\mathbb{Q})\rightarrow H^*
(\widetilde{T}_i,\mathbb{Q})$$
are surjective for all $i$. It thus follows that the cohomology of $Y$ is generated (as an algebra)
by
the pull-back of the cohomology of $\widetilde{T\times T}_{x_1,\ldots,x_N}$ and the classes
$e_i$ of the exceptional divisors over $\widetilde{T}_i$.

We can thus apply Theorem \ref{theoprojbun} and conclude
that  for any Hodge structure on the rational cohomology algebra
$H^*(X,\mathbb{Q})$, the sub-algebra $\pi^*H^*(Y,\mathbb{Q})$ is a sub-Hodge structure, and
furthermore, for the induced rational Hodge structure on
$H^*(Y,\mathbb{Q})$, the Chern class $c_2(L\oplus -L)=-\lambda^2$ is a Hodge class.

We now briefly recall the analysis made in \cite{voisin}: We proved  that
for any rational Hodge structure on $H^*(Y,\mathbb{Q})$, the classes $e_i$ are Hodge classes.
Studying the cup-product maps (which must be morphisms of Hodge structures of bidegree $(1,1)$)
$$e_i\cup :H^1(Y,\mathbb{Q})\rightarrow H^3(Y,\mathbb{Q}),$$
where $H^1(Y,\mathbb{Q})\cong H^1(T\times T,\mathbb{Q})$, we
then  concluded :
\begin{enumerate}\item \label{1} The induced
(weight $1$)  Hodge
structure on $H^1(T\times T,\mathbb{Q})$ is induced by a Hodge structure on
$H^1(T,\mathbb{Q})$, that is,  is the direct sum of two copies of a rational Hodge structure
on $H^1(T,\mathbb{Q})$.

\item \label{2} Furthermore this Hodge structure must admit the automorphism
 $^t\phi$. (Observe that $H^1(T,\mathbb{Q})\cong \Gamma^*_\mathbb{Q}$, so that $^t\phi$
 acts on it.)
 \end{enumerate}

The proof is now finished because we know  that the class $\lambda^2$ must be a Hodge class
in $H^4(Y,\mathbb{Q})$, and we made the assumption that it belongs to
the sub-Hodge structure
$\tau^*H^4(T\times T,\mathbb{Q})=\bigwedge^4H^1(Y,\mathbb{Q})
\subset H^4(Y,\mathbb{Q})$. Thus it must be a non trivial Hodge classes in
$H^4(T\times T,\mathbb{Q})$, for the Hodge structure induced by the Hodge structure above on
$H^1(T,\mathbb{Q})$. However, as in \cite{voisin} (see
\cite{voisinjdg} for more details of such computations), using assumption (\ref{condition*}),
an easy irreducibility argument for the action
of $^t\phi$ on certain natural direct summands of
$H^*(T,\mathbb{Q})\otimes H^*(T,\mathbb{Q})$  shows that there is no non zero
Hodge classes in $H^4(T\times T,\mathbb{Q})$ for any Hodge structure satisfying the conclusions
\ref{1}, \ref{2} above.

\cqfd

\subsection{Constraints coming from polarizations \label{sectionpol}}
We have made  essentially no use of the polarization in the previous sections. Furthermore, we used a
lot odd dimensional cohomology, for the following reason: a Hodge structure of odd weight
always gives  a non trivial information. This is not the case
for Hodge structures of even weight
$2k$. In this case, we can always consider the trivial Hodge structure, for which everything is of type
$(k,k)$. Furthermore, if we have a cohomology algebra $M$ with trivial odd degree part, then
we can always put the trivial Hodge structure on each term
$M^{2k}$, and this will give us a cohomology algebra with Hodge structure.

In the presence of polarization, and assuming the dimension
of $M$ is divisible by $4$,  we meet now the following restriction given by the Hodge index theorem
\ref{rieho}. Let us define the signature $\tau(M)$ of
a cohomology algebra  $M$ of dimension $4m$ as the signature on the intersection form
on $M^{2m}$. This signature is only defined up to sign if we did not fix the isomorphism
$$M^{4m}\cong \mathbb{R}.$$

\begin{theo}\label{rieho} Let $M$ be a cohomology algebra of dimension $4m$
which is endowed with a trivial Hodge structure. Then, if this Hodge structure can be polarized,
 we have
$$\tau(M)=\pm \sum_i(-1)^i rk\,M^{2i}.$$
\end{theo}
{\bf Proof.}
Indeed, if we look at the Riemann second bilinear relations (\ref{secondriemann}), they give us
the signs on the intersection form on the real part of
the $(p,q)\oplus (q,p)$ part of the pieces of the Lefschetz decomposition
on the middle piece
$M^{2n}$. (Here, all the signs may be reversed by changing the orientation.)
As a result one gets the Hodge index formula
\cite{voisinbook}, 6.3.2, which gives the signature
$\tau$ as an alternate sum of $m^{p,q}$ numbers, $m^{p,q}:=rk\,M^{p,q}$.
In our case, the Hodge decomposition is trivial,
so  we  get the formula
\begin{eqnarray}\label{formulasign}\tau=\pm\sum_i(-1)^im_{2i},\,m^{2i}:=rk\,M^{2i}.
\end{eqnarray}
\cqfd

Let us now use this condition, combined with Deligne's Lemma
\ref{dele},  to construct a symplectic compact manifold $X$ such that
$M=H^*(X,\mathbb{Q})$ carries
a rational Hodge structure, and such that adequate symplectic classes $\omega$ on $X$
satisfy the Lefschetz property, but such that $M$ does not carry any {\it polarized} Hodge structure.
\begin{ex}
We consider a $K3$ surface $S$, and consider a basis
$a_1,\,\ldots, a_{22}$ of $H^2(S,\mathbb{Q})$ consisting of symplectic classes; more precisely,
we assume that $a_i$ is the cohomology class of a symplectic form
$\alpha_i$ close to a given symplectic form
$\alpha$, so that $\sum_i\alpha_i$ is again a symplectic form.

For  an adequate integer coefficient $l>>0$, one knows
by  \cite{gromov}  that
one can construct embeddings $\phi_i:S\rightarrow \mathbb{P}^N$ such that
$\phi_i^*\Omega=l\alpha_i$, where $\Omega$
is the Fubini-Study K{\"a}hler form  on $\mathbb{P}^N$. Then
$\sum_i\phi_i^*\Omega=l\sum_i\alpha_i$ is again a symplectic form on $S$.

We consider now the  embedding
$\psi:=(\phi_1,\ldots,\phi_{22})$ of $S$ into $(\mathbb{P}^N)^{22}$.
With the above choice of $\phi_i$'s, this provides a symplectic submanifold of
$(\mathbb{P}^N)^{22}$ endowed with the product K{\"a}hler form.

The symplectic manifold we will consider will be
the symplectic blow-up of $(\mathbb{P}^N)^{22}$ along $\psi(S)$.
This is a symplectic manifold with symplectic class given by
$\tau^*(\sum_i pr_i^*\omega)-\epsilon e$, where $\tau:X\rightarrow(\mathbb{P}^N)^{22}$ is the
blowing-up map and $e$ is the cohomology class of the exceptional divisor.
\end{ex}
Let us show:
\begin{theo} \label{theopol}
Any Hodge structure on $H^*(X,\mathbb{Q})$ is trivial (that is
$H^{2i}(X,\mathbb{Q})$ is completely of type $(i,i)$). Furthermore the trivial Hodge structure
cannot be polarized.
\end{theo}
{\bf Proof.} We use lemma \ref{dele} to show that for any Hodge structure
on $H^*(X,\mathbb{Q})$, the cohomology $H^2(X,\mathbb{Q})$ must be completely of type
$(1,1)$. Indeed, the cohomology of $X$ in degree $2$ is generated by the
$h_i:=\tau^*pr_i^*h$ and by $e$.
We observe that for each $i=1,\ldots,22$,
the line generated by  $h_i$ is an irreducible component of the set
$$Z\subset H^2(X,\mathbb{C}),\,Z=\{\alpha\in H^2(X,\mathbb{C}),\,\alpha^{N+1}=0\}.$$
Thus Lemma \ref{dele} implies that this line is generated by a Hodge class, as it is defined over
$\mathbb{Q}$. Finally, the last generator $e$ must also be of type $(1,1)$, because
the $(2,0)\oplus (0,2)$ part has even rank, and we just proved above that it has rank $\leq 1$.

Next we observe that the cohomology of our variety $X$ is generated in degree $2$.
This follows from  the computation of the cohomology of a blown-up variety (cf \cite{voisinbook}, 7.3.3
in the complex case,
the symplectic case is computed in the same way), and from the
following  facts:

1) The variety $(\mathbb{P}^N)^{22}$ satisfies the property that its cohomology is generated in degree
$2$.

2) The restriction map on  cohomology
$$\psi^*:H^*((\mathbb{P}^N)^{22},\mathbb{Q})\rightarrow H^*(S,\mathbb{Q})$$
is surjective.

It follows that if $j:E\hookrightarrow X$ is the inclusion, and $\tau':E\rightarrow S$
is the restriction of $\tau$, then
for $\alpha=\psi^*\beta\in H^*(S,\mathbb{Q})$, we have
$$j_*\tau^{'*}\alpha=e\cup\tau^*\beta\,\,{\rm in}\,\,H^*(X,\mathbb{Q}),$$
and this implies that the cohomology of $X$ is also generated in degree $2$.

As the Hodge structure on $H^2(X,\mathbb{Q})$ is trivial, the Hodge structure on all
cohomology groups
$H^{2i}(X,\mathbb{Q})$ are trivial too. Furthermore, $X$ has no odd dimensional cohomology.

It remains to see why the trivial Hodge structure cannot be polarized. For this we apply
Theorem \ref{rieho}, which tells us that in this case we should have
\begin{eqnarray}\label{formulasign}\tau(X)=\pm\sum_i(-1)^ib_{2i}.
\end{eqnarray}
This formula must give   the signature (defined up to sign)
of  the intersection form on the middle part
of any cohomology algebra endowed with a trivial Hodge structure which can be polarized.
We can now easily construct a projective variety with the same Betti numbers
and which has the property that the Hodge structures on its cohomology are trivial.
Namely we start from $S'=\mathbb{P}^2$ blown-up at $21$ points,
choose a basis $b_1,\ldots,\,b_{22}$ of $H^2(X',\mathbb{Q})$ consisting of
Chern classes of very ample line bundles $L_i$ and imbed holomorphically
$S'$ in $(\mathbb{P}^N)^{22}$, using the morphisms $\phi'_i$ given by
the line bundles $L_i$.
Then we consider the projective  variety $X'$ defined by blowing-up
$S'$ in $(\mathbb{P}^N)^{22}$.

As $X'$ is projective, its cohomology carries a polarized Hodge structure. As the blown-up
surface $X'$ has a trivial Hodge structure, the Hodge structure on $H^*(X',\mathbb{Q})$
is trivial. Note also that, as $S'$ has only trivial Hodge structures, so does $X'$,
and thus
$X'$ has the same Hodge numbers as $X$.
 Thus the signature of $X'$ is given by formula
 (\ref{formulasign}).

 Hence, to conclude that the trivial Hodge structure on $H^*(X,\mathbb{Q})$ cannot be polarized,
 it suffices to show that the absolute value of the signature of $X$ is different from that of
 $X'$, hence does not satisfy formula (\ref{formulasign}).

 This is quite easy, because the middle cohomology
 $H^{22N}(X,\mathbb{Q})$ is a direct   sum
 $$H^{22N}(\mathbb{P}^{22N},\mathbb{Q})
 \bigoplus e^{11N}\mathbb{Q}\bigoplus j_*(e^{11N-3}\tau^{'*}H^4(S,\mathbb{Q}))\bigoplus
 j_*(e^{11N-2}\tau^{'*}H^2(S,\mathbb{Q})),$$
 where the first term is orthogonal to the three other ones, the second and the third are
 isotropic and dual, and orthogonal to the last one, and the intersection form of $X$ restricted
 to
 the term $j_*(e^{11N-2}\tau^{'*}H^2(S,\mathbb{Q})$ is equal to the intersection form on
 $H^2(S,\mathbb{Q})$, with opposite sign.

We do the same computation with $X'$ and we conclude that the difference
$\tau(X)-\tau(X')$ is equal to
$\tau(S')-\tau(S)$, hence is non zero.

The argument is not quite complete, as we also have to show
that we do not have $\tau(X)=-\tau(X')$.   This is easily  checked for large $N$.

\cqfd

\section{Further restrictions}
Up to now, we have been studying the topological constraints on compact K\"ahler manifolds
via the  constraints imposed to the cohomology algebra by the existence  of
a (polarized) Hodge structure. In this final section, we want to show that there are
in fact other constraints on the cohomology algebra.
\begin{theo} There exists a compact symplectic manifold whose rational  cohomology algebra carries a
rational polarizable Hodge structure, but is not isomorphic to the rational
cohomology algebra of any compact K\"ahler manifold.
\end{theo}
In fact our  example will  even be  a manifold which does not have the rational cohomology
algebra of a compact K\"ahler manifold,  but satisfies the following properties:
\begin{enumerate}
\item The cohomology algebra of $X$ admits a polarized rational Hodge structure.
\item $X$ has  the real cohomology algebra of
a compact K{\"a}hler manifold.

\end{enumerate}
\begin{ex}
We consider a $3$-dimensional torus $T$, which admits complex multiplication by a number field
$K$, with $[K:\mathbb{Q}]=6$. This means  that
$K$ acts on $T$, hence on $H^1(T,\mathbb{Q})$, and this  makes the space
$H^1(T,\mathbb{Q})$  a $1$-dimensional $K$-vector space.

Let $k_i,\,i=1,\ldots, 6$ be a basis of $K$ over $\mathbb{Q}$ and for each
$i$, let
$$\gamma_i\in Hom\,(H^1(T,\mathbb{Q}),H^1(T,\mathbb{Q}))\subset
H^6(T\times T,\mathbb{Q})$$
 be  given by the action of $k_i$
on $H^1(T,\mathbb{Q})$. Let also $f_1,\,f_2\in H^6(T\times T,\mathbb{Q})$
be the respective classes of the fibers of the projections
$pr_1:T\times T\rightarrow T$, $pr_2:T\times T\rightarrow T$.

Note that we have
\begin{eqnarray}\label{ouf} H^6(T\times T,\mathbb{Q})\otimes H^6(T\times T,\mathbb{Q})\subset
H^{12}(T\times T\times T\times T,\mathbb{Q})\end{eqnarray}
by K\"unneth decomposition, and
\begin{eqnarray}\label{oufouf}H^6(T\times T,\mathbb{Q})\otimes H^6(T\times T,\mathbb{Q})\cong Hom\,
(H^6(T\times T,\mathbb{Q}), H^6(T\times T,\mathbb{Q}))\end{eqnarray}
by Poincar{\'e} duality on $T\times T$.

Let $V\subset H^6(T\times T,\mathbb{Q})$ be the subspace generated by
$f_1,f_2,\gamma_1,\ldots,\gamma_6$. Observing that the intersection pairing is non degenerate on
$V$, let $$P:H^6(T\times T,\mathbb{Q})\rightarrow  H^6(T\times T,\mathbb{Q})$$
denote the orthogonal projector onto $V$.

By (\ref{ouf}) et (\ref{oufouf})
$P$ can be seen as a rational cohomology class of degree $12$ on $Y:=T^4$.
Furthermore, $P$ is a Hodge class, because it corresponds to an  endomorphism of  Hodge structure
(an orthogonal projection onto a sub-Hodge structure) acting
on $H^*(T\times T,\mathbb{Q})$.

Let $G$ be a complex vector bundle of rank $r$ on $Y$ such that
$c_i(G)=0$ for $i\not=0,\,6$, and $c_6(G)$ is a non zero multiple of $P$, say
$c_6(G)=mP$.

Similarly the classes $e,\, f$ of the fibers  of the projections $pr_1$, resp. $pr_2$
from $Y=T^4=T^2\times T^2$ to $T^2$ are degree $12$ Hodge classes on $Y$.
Let $E,\,F$ be complex vector bundles of rank $s$ on $Y$ with
 the property $c_i(E)=c_i(F)=0,\,i\not=0,\,6$, $c_6(E)=m' e$, $c_6(F)=m'f$ for some non zero integer
 $m'$.

We define $X$ to be the fibered product over $Y$ of
the projective bundles $\mathbb{P}(E),\,\mathbb{P}(F),\,\mathbb{P}(G)$:
$$X:=\mathbb{P}(E)\times_Y \mathbb{P}(F)\times_Y\mathbb{P}(G).$$
\end{ex}

Let us prove:
\begin{theo}\label{theofinal}
\begin{enumerate}
\item\label{item115mars}
 The cohomology algebra $H^*(X,\mathbb{Q})$ is endowed with a natural Hodge structure.
 \item \label{item215mars} The real cohomology algebra $H^*(X,\mathbb{R})$
 does not  depend on the choice of number field
$K$ satisfying the condition that $K\otimes\mathbb{R}\cong \mathbb{C}^3$.
\item \label{item315mars} For certain choices of $K$, and adequate choices
of $r,\,s, \,m,\,m'$, the variety $X$ is projective
(in particular K{\"a}hler). With the same $m,\,m',\,r,\,s$, for a ``generic'' choice of $K$, $X$
does not have the rational cohomology algebra of a K{\"a}hler manifold.
\end{enumerate}
\end{theo}

  {\bf Proof.}
  Statement \ref{item115mars} follows from the explicit computation of the cohomology algebra of
  $X$, which is the fibered product of the projective bundles
  $\mathbb{P}(E),\,\mathbb{P}(F),\,\mathbb{P}(G)$ over $Y=T^4$.
  This cohomology algebra is then generated over
  $H^*(T^4,\mathbb{Q})$ by
  $h_{E}:=c_1(\mathcal{L}_{\mathbb{P}(E)})$, where $\mathcal{L}_{\mathbb{P}(E)}$ is the
  dual of the relative Hopf line bundle, and similarly $h_F$ and $h_G$,
   with the relations
  $$h_{E}^s=-\pi^*c_6(E)h^{s-6}_{E},\,\,\,h_F^{s}=-\pi^*c_6(F)h^{s-6}_{F},\,\,\,
  h_G^r=-\pi^*c_6(G)h^{r-6}_{G},$$
where  $\pi:X\rightarrow Y$
is the structural map. This presentation is due to the fact that
$c_i(E)=0,\,i\not=0,\,6$ and similarly for $F$ and $G$.

As we know that $c_6(E),\,c_6(F), \,c_6(G)$ are Hodge classes on $Y$, there is a Hodge structure on
$H^*(X,\mathbb{Q})$, inducing the Hodge structure on $H^*(Y,\mathbb{Q})$
and  obtained by declaring
$h_{E},\,h_F$ and $h_G$ to be of type $(1,1)$.

For the proof of \ref{item215mars},  observe that the real cohomology algebra
$H^*(X,\mathbb{R})$ depends only on the
representation of the algebra $K\otimes_\mathbb{Q}{\mathbb{R}}$
on $H^1(T,\mathbb{R})$. Thus  we conclude that the $\mathbb{R}$-algebra
$H^*(X,\mathbb{R})$ does not depend on the choice of field $K$ satisfying the
condition that $K\otimes_\mathbb{Q}\mathbb{R}=\mathbb{C}^3$.

It remains to prove \ref{item315mars}.

First of all, we note that we can choose $K$ in such a way that
$T$ and thus $Y=T^4$ are abelian
varieties. The classes $e,\,f, f_1,\,f_2,\,\gamma_i$ are all classes of algebraic cycles.
Here we use the fact that $\gamma_i$ are algebraic, which is a general fact (see \cite{murre}):
the K\"unneth component of type $(n-1,1)$ of a codimension
$n$  algebraic cycle class
in a product of smooth projective varieties $Z,Z',\,dim\,Z=n$, is again
an algebraic cycle class.

As $Y$ is projective, we know that the Chern character
$$ch: K_{0,alg}\otimes\mathbb{Q}\rightarrow CH(Y)\otimes\mathbb{Q}$$
is an isomorphism, and this implies that for adequate choices of ranks
$r,s$ and integers $m,\,m'$
there exist {\it algebraic} vector bundles $E,\,F$ of rank $s$, $G$ of rank
$r$ on $Y$
satisfying the conditions :
$$c_i(E)=c_i(F)=c_i(G)=0,\,i\not=6,\,c_6(E)=m'e,\,c_6(F)=m'f,\,c_6(G)=m P.$$
The corresponding manifold
$$X=\mathbb{P}(E)\times_Y\mathbb{P}(F)\times_Y\mathbb{P}(G)$$
is then projective, which proves the first statement.

Let us now show that for a number field $K$ satisfying the condition
(\ref{condition*}), that is, its Galois group is as large as possible, the
rational cohomology algebra of the symplectic  manifold
$X$ (independently of the integers $r,\,s,\,m,\,m'$) is not isomorphic to the rational cohomology
algebra of a compact K\"ahler manifold. The key point here, as in \cite{voisin} or
in Theorem \ref{theotorus} above,  is that in this case,
the presence of such an algebra acting by isogenies on the complex torus $Y$ prevents it to be an abelian
variety. In turn, this will prevent the Hodge classes used below
to come from Chern classes of reflexive analytic coherent sheaves, as in \cite{voisinimrn}.

So assume to the contrary that there exists a compact K\"ahler manifold $Z$
which has its cohomology algebra isomorphic to that of $X$. Let
$\pi':Z\rightarrow Alb\,Z$ be the Albanese map of $Z$. Topologically, $\pi'$ induces
the isomorphism
$$\pi^{'*}:H^1(Alb\,Z,\mathbb{Q})\cong H^1(Z,\mathbb{Q})$$
and the map induced by cup-product on $Z$:
$$\pi^{'*}:H^l(Alb\,Z,\mathbb{Q})\cong \bigwedge^lH^1(Z,\mathbb{Q})\rightarrow H^l(Z,\mathbb{Q}).$$
As the cohomology algebra of $X$ and $Z$ are isomorphic, we conclude that
the map $\pi^{'*}$ is injective in top degree for $Z$, as it is the case for $X$, and this implies that
$$\pi':Z\rightarrow Alb\,Z$$
is surjective.

Observe also for future use that, via the isomorphism of cohomology algebras
$$H^*(X,\mathbb{Q})\cong H^*(Z,\mathbb{Q})$$
the morphism $\pi^{'*}$ identifies to $\pi^*:H^*(Y,\mathbb{Q})\rightarrow H^*(X,\mathbb{Q})$.

Next, we have the following
lemma:
\begin{lemm}\label{delafin}
 The torus $Alb\,Z$ must be isogenous to
a product $T^{2}_1\times T_2^2$ where $T_1$
and $T_2$ are   three-dimensional tori on which the field $K$ acts by isogenies.
\end{lemm}
We postpone the proof of this lemma and conclude the proof of Theorem \ref{theofinal}
as follows: As the torus $Alb\,Z$ is of the form $T_1^2\times T_2^2$, where
the $T_i$'s are  $3$-dimensional tori with an action of $K$, and as $K$
satisfies condition (\ref{condition*}), one proves as in
\cite{voisinjdg} that
$$Hdg^2(Alb\,Z)=0,\,Hdg^4(Alb\,Z)=0.$$
This torus thus satisfies the condition of the Appendix of \cite{voisinimrn}, and we thus conclude
as in \cite{voisinimrn} that
any reflexive analytic  coherent sheaf $\mathcal{F}$ on $Alb\,Z$
 is a vector bundle with trivial Chern classes.
In particular, by the Riemann-Roch formula for
complex  vector bundles
on compact complex manifolds (see \cite{toledotong}),
any reflexive coherent sheaf on $Alb\,Z$ satisfies
\begin{eqnarray}\label{annutorecoh}\chi(Alb\,Z,\mathcal{F})=0.
\end{eqnarray}
Next we observe that the morphism
$\pi': Z\rightarrow Alb\,Z$ is projective, that is, there exist
 line bundles on $Z$ which are relatively ample w.r.t. $\pi'$.
Indeed, this follows from iterated applications of
(the proof of) Theorem \ref{theoprojbun}, which show that
for any Hodge structure on the
rational cohomology algebra $H^*(Z,\mathbb{Q})$, the
sub-algebra $Im\,\pi^{'*}$ (which identifies as mentioned above to $Im\,\pi^*$)
is a sub-Hodge structure, and that
$H^2(Z,\mathbb{Q})$ is generated by $\pi^{'*}H^2(Alb\,Z,\mathbb{Q})$ and by three Hodge classes
of degree $2$ (corresponding to $h_E,\,h_F,\,h_G$).
It follows that for any  fiber $Z_a$ of $\pi^{'*}$,  the image of
the map $H^2(Z,\mathbb{Q})\rightarrow H^2(Z_a,\mathbb{Q})$ is of type $(1,1)$.
On the other hand, this image contains a K\"ahler class on $Z_a$.
We then easily conclude that some rational combinations of the classes
$h_E,\,h_F,\,h_G$ (transported to $Z$) restricts to a K\"ahler class on any fiber
$Z_a$, which implies by Kodaira theorem that this rational combination
is ample on the fibers of $\pi'$.

The contradiction now comes from the following:
let $\mathcal{L}$ be a relatively ample line bundle on $Z$, whose first Chern class
$c_1(\mathcal{L})$
is a rational combination of $h_E,\,h_F,\,h_G$. Observe that we can assume that the top
self-intersection
$c_1(\mathcal{L})^N,\,N=dim\,Z=12+r-1+2(s-1)$ is not equal to $0$. Indeed, this follows from
the relations
$$h_G^r=c_6(G)h_G^{r-6},\,\,\,h_G^{r+11}=c_6(G)h_G^{r+5}=
c_6(G)^2h_G^{r-1},$$
$$h_G^{r+11}h_E^{s-1}h_F^{s-1}=c_6(G)^2h_E^{s-1}h_F^{s-1}h_G^{r-1}\not=0,$$
where the last equality follows from $c_6(G)^2=m^2P^2\not=0$.

As $c_1(\mathcal{L})^N\not=0$, the Hilbert  polynomial
$P_\mathcal{L}$ defined by
$P_\mathcal{L}(n)=\chi(Z,\mathcal{L}^{\otimes n})$ is not identically equal to $0$.
As $\mathcal{L}$ is relatively ample, we have
$$R^i\pi'_*\mathcal{L}^{\otimes n}=0$$ for large $n$ and $i>0$, and thus

$$\chi(Z,\mathcal{L}^{\otimes n})=\chi(Alb\,Z,R^0\pi^{'*}\mathcal{L}^{\otimes n}).$$
But the coherent  sheaf $\mathcal{F}_n=R^0\pi'_*\mathcal{L}^{\otimes n}$ is reflexive on
$Alb\,Z$, because the analysis of $H^2(Z)$ shows that no divisor
of $Z$ can be contracted to a codimension $2$ analytic subset of
$Alb\,Z$. Hence by (\ref{annutorecoh}),  $\mathcal{F}_n$ satisfies
$\chi(Alb\,Z,\mathcal{F}_n)=0$, which is a contradiction.

\cqfd
{\bf Proof of Lemma \ref{delafin}.} We use first of all
in an iterated way Theorem \ref{theoprojbun}, which implies that
any Hodge structure on $H^*(X,\mathbb{Q})\cong H^*(Z,\mathbb{Q})$
 comes from a Hodge structure
on $H^*(Y,\mathbb{Q})$ for which the classes $c_6(E),\,c_6(F),\,c_6(G)$ must be
Hodge classes. As we have already seen,
$H^*(Y,\mathbb{Q})$ identifies to $H^*(Alb\,Z,\mathbb{Q})$ under the isomorphism
$$H^*(X,\mathbb{Q})\cong H^*(Z,\mathbb{Q}).$$
Thus the complex torus $Alb\,Z$ satisfies the property
that there is an isomorphism of exterior algebras
$$ H^*(T^4,\mathbb{Q})\cong H^*(Alb\,Z,\mathbb{Q})$$
sending the classes $e,\,f$ and $P$ to Hodge classes.

Recall that $Y\cong T^2\times T^2$ and that $e$ and $f$ are the fibers of the two projections
on $T^2$.
It follows that
$$Ker \,(e\wedge :H^1(Y,\mathbb{Q})\rightarrow H^1(Y,\mathbb{Q})),$$
which is equal to $pr_1^*H^1(T^2,\mathbb{Q})$ is sent to a sub-Hodge structure
of $H^1(Alb\,Z,\mathbb{Q})$. Similarly,
$Ker \,(f\wedge :H^1(Y,\mathbb{Q})\rightarrow H^1(Y,\mathbb{Q}))$,
which is equal to $pr_2^*H^1(T^2,\mathbb{Q})$ is sent to a sub-Hodge structure
of $H^1(Alb\,Z,\mathbb{Q})$. Thus we conclude that
$Alb\,Z$ is isogenous to a direct sum
$T'\oplus T''$ of two tori, in such a way that
the isomorphism
$$H^1(T^2\times T^2,\mathbb{Q})\rightarrow H^1(Alb\,Z,\mathbb{Q})$$
sends $pr_1^*H^1(T^2,\mathbb{Q})$ to $H^1(T',\mathbb{Q})$ and
$pr_1^*H^1(T^2,\mathbb{Q})$ to $H^1(T'',\mathbb{Q})$.

We now consider the image of the class $P$. This is now a Hodge class
in $H^*(T'\times T'',\mathbb{Q})$ which lies in $Hom\,(H^*(T',\mathbb{Q}),H^*(T'',\mathbb{Q}))$.
Thus its image in $ H^*(T'',\mathbb{Q})$ and the orthogonal of its kernel
in $ H^*(T',\mathbb{Q})$, which both identify to
$V\subset H^*(T^2,\mathbb{Q})$, are sub-Hodge structure of
$H^*(T'',\mathbb{Q})$ and $H^*(T',\mathbb{Q})$ respectively. Thus it suffices to prove that
if a $6$-dimensional torus $T'$ admits an
isomorphism of cohomology algebras
$$H^*(T\times T,\mathbb{Q})\cong H^*(T',\mathbb{Q})$$
sending $V$ to a sub-Hodge structure of $H^*(T',\mathbb{Q})$, then
$T'$ is of the form $T_1^2$, where $T_1$ admits an action of $K$ by isogenies.

We first observe
that the complex lines generated by $f_1$ and $f_2$ are irreducible components of
the set of reducible elements in
$V_\mathbb{C}\subset \bigwedge^6H^1(T',\mathbb{C})=H^6(T',\mathbb{C})$.
As these two lines are defined over $\mathbb{Q}$,
they must be generated by a Hodge class on $T'$, by lemma \ref{dele}.

We use now  the classes $f_1,\,f_2$ to show that $T'$ is isogenous to a product
$$T'\cong T'_1\oplus T'_2$$
 of two $3$-dimensional complex tori. Indeed, we recover the sub-Hodge structures $H^1(T'_i,\mathbb{Q})
 \subset H^1(T',\mathbb{Q})$ as the kernel of the morphism of
 Hodge structures
$$\cup f_i: H^1(T',\mathbb{Q})\rightarrow H^7(T',\mathbb{Q}).$$

Having this decomposition, we look at the image of $V$ in
$Hom\,(H^1(T'_1,\mathbb{Q}), H^1(T'_2,\mathbb{Q}))$
via the natural projection (given by K\"unneth decomposition and Poincar{\'e} duality)
$$H^6(T'_1\times T'_2,\mathbb{Q})\rightarrow
Hom\,(H^1(T'_1,\mathbb{Q}), H^1(T'_2,\mathbb{Q})).$$

The image $\overline{V}$ of $V$ under this map is a sub-Hodge structure  of rank $6$, and
the subalgebra $L$ of $End\,H^1(T'_1,\mathbb{Q})$ generated by the $\gamma_i^{-1}\circ\gamma_j$
is a sub-Hodge structure of $End\,H^1(T'_1,\mathbb{Q})$ which is isomorphic to $K$
as a $\mathbb{Q}$-algebra.
As $K\otimes\mathbb{C}$ has no nilpotent element, it follows that $L_\mathbb{C}$
has no $(-1,1)$-part in its Hodge decomposition, because the
$(-1,1)$-part of the Hodge structure on $End\,H^1(T'_1,\mathbb{Q})$ is equal to
$Hom\,(H^{0,1}(T'_1),H^{1,0}(T'_1))$ and this is nilpotent.

Thus $L_\mathbb{C}$ is purely of type $(0,0)$, and $L$ consists of endomorphisms of $T'_1$.
As $L$ is isomorphic to $K$, this shows that $T'_1$ has complex multiplication
by $K$. Finally, as $L_\mathbb{C}$ is of type $(0,0)$, the same is true of
$\overline{V}$, which implies that $T'_1$ and $T'_2$
are isogenous. This concludes the proof of the Lemma.
\cqfd

\end{document}